\documentclass[11pt]{article}
\usepackage[margin=1in]{geometry}
\usepackage{amsmath,amssymb,amsthm,mathtools,microtype}
\usepackage[colorlinks=true,linkcolor=blue,citecolor=blue,urlcolor=blue]{hyperref}
\newtheorem{theorem}{Theorem}[section]
\newtheorem{lemma}[theorem]{Lemma}
\newtheorem{proposition}[theorem]{Proposition}
\newtheorem{corollary}[theorem]{Corollary}
\theoremstyle{remark}

\newcommand{\R}{\mathbb{R}}
\newcommand{\Z}{\mathbb{Z}}
\newcommand{\N}{\mathbb{N}}

\newcommand{\E}{\mathbb{E}}

\newcommand{\e}{\mathrm e}
\newcommand{\1}{\mathbf 1}

\newcommand{\ga}{\gamma}

\newcommand{\eps}{\varepsilon}

\renewcommand{\Re}{\operatorname{Re}}

\title{On the van der Corput property for cubes}
\author{%
	Steve Fan\textsuperscript{1}, Andrew Lott\textsuperscript{1},
	Ákos Magyar\textsuperscript{1}, Máté Matolcsi\textsuperscript{2,3},
	and Imre Z.~Ruzsa\textsuperscript{2}\\[0.55em]
	\small\textsuperscript{1}\,\textit{University of Georgia}\\
	\small\textsuperscript{2}\,\textit{HUN-REN Alfréd Rényi Institute of Mathematics}\\
	\small\textsuperscript{3}\,\textit{Budapest University of Technology and Economics}%
}
\date{}
\begin{document}
	\maketitle
	\begin{abstract}
		\noindent
		We prove a power-saving form of the van der Corput property for cubes. As a consequence, we obtain a power-saving bound $|A|\ll_\varepsilon N^{59/60+\varepsilon}$ for sets $A\subseteq\{1,\ldots,N\}$ with no nonzero cube differences.
	\end{abstract}
	
	\section{Introduction}
	Answering a question of Lovász, Furstenberg~\cite{Fur} and Sárközy~\cite{Sar}
	independently proved that if $A\subseteq[N]:=\{1,\ldots,N\}$ has no nonzero
	square difference, then $|A|=o(N)$. Sárközy's proof gave
	\[
	|A|\ll N(\log\log N)^{2/3}(\log N)^{-1/3}.
	\]
	Pintz, Steiger, and
	Szemerédi~\cite[Theorem 1]{PSS} used a more elaborate version of Sárközy's density increment strategy to prove
	\[
	|A|\ll N(\log N)^{-\frac{1}{12}\log\log\log\log N}.
	\]
	Bloom and Maynard~\cite{BM} later improved this bound to
	\[
	|A|\ll N(\log N)^{-c\log\log\log N}
	\]
	for some absolute constant $c>0$. Recently, Green and Sawhney~\cite{GS} combined Sárközy's density increment strategy with their arithmetic level-$d$ inequality to show that, for some absolute $c>0$,
	\[
	|A|\ll N\exp(-c\sqrt{\log N}).
	\]
	Fan and Lott~\cite{FL} subsequently used the van der Corput property for squares to improve this to
	\begin{equation}\label{eq:fan-lott-squares}
		|A|\ll N\exp\!\bigl(-c\sqrt{\log N\log\log N}\bigr)
	\end{equation}
	for some absolute constant $c>0$.
	
	A longstanding question is whether there is an absolute constant $c>0$ such that
	every $A\subseteq[N]$ with no nonzero square difference satisfies the power-saving bound $|A|\ll N^{1-c}$. Very recently, OpenAI~\cite{OAIsq} released a preprint claiming such a bound.

	The lower bound of Fan and Lott~\cite[Theorem~1.3]{FL} shows that the van der Corput method cannot yield a power saving for sets with no nonzero square
	difference. For cubes, however, we obtain a power saving by exploiting the oddness of $x\mapsto x^3$, which makes the complete cubic exponential sums real.
	
	Write $\e(t)=\exp(2\pi i t)$. Let $F\subseteq\Z$, let $N$ be a positive integer, and let
	$\delta>0$. Consider a real trigonometric polynomial
	\begin{equation}\label{eq:intro-vdc-polynomial}
		P(\theta)=a_0+\sum_{n=1}^N a_n\cos(2\pi n\theta),
		\qquad a_0,a_1,\ldots,a_N\in\R,
	\end{equation}
	such that
	\begin{alignat*}{2}
		\textnormal{(i)}\quad   &P(\theta)\ge0 &&\qquad\text{for every }\theta\in\R/\Z,\\
		\textnormal{(ii)}\quad  &P(0)=1,\\
		\textnormal{(iii)}\quad &a_n=0 &&\qquad\text{if }1\le n\le N\text{ and }n\notin F\cup(-F),\\
		\textnormal{(iv)}\quad  &0\le a_0\le\delta.
	\end{alignat*}
	We call a polynomial of the form~\eqref{eq:intro-vdc-polynomial}
	satisfying \textnormal{(i)--(iii)} a \emph{van der Corput witness} for
	$F$ of degree at most $N$. Its nonconstant coefficients may have either
	sign. In the trigonometric characterization of Kamae and Mendès
	France~\cite{KMF} and Ruzsa~\cite{Ruzsa}, $F$ has the \emph{van der
		Corput property} if for every $\delta>0$ there exist an
	integer $N$ and a witness for $F$ satisfying \textnormal{(iv)}.

	Suppose that $A\subseteq[N]$ and $(A-A)\cap F\subseteq\{0\}$. For any $P$
	satisfying \textnormal{(i)--(iv)}, Fourier orthogonality gives
	\begin{equation}\label{eq:intro-vdc-bound}
		|A|\le (2N+1)a_0\le(2N+1)\delta.
	\end{equation}
	We give the argument in the proof of Corollary~\ref{cor:cube-differences}.

	For the set of shifted primes $F=\{p-1:p\ \text{prime}\}$,
	Green~\cite[Theorems~1.1 and~1.2]{GreenShifted} constructed a witness
	of degree at most $N$ with $a_0\ll N^{-c}$ for some absolute $c>0$.
	It follows that $|A|\ll N^{1-c}$ whenever $A\subseteq[N]$ contains
	no two elements differing by one less than a prime.

	Several other sets of forbidden differences admit power savings by the
	van der Corput method. Fan and Lott~\cite{FLTwo} proved that
	if $A-A$ contains no nonzero sum of two squares, then
	$|A|\ll_\varepsilon N^{7/8+\varepsilon}$ for every $\varepsilon>0$. Younis~\cite{Younis} proved that $|A|\ll N^{1-c}$
	for some absolute $c>0$ if $A-A$ contains no integer of the form
	$n\lfloor 2^{1/3}n\rfloor$ with $n\in\N=\{1,2,\ldots\}$.
	Kalogirou, Lott, Magyar, and Singha Roy~\cite{KLMS} proved that for
	each odd $k\ge3$ and even $s\ge2$ there is a $c_{k,s}>0$ such that
	$|A|\ll_{k,s}N^{1-c_{k,s}}$ if $A-A$ contains no sum of $s$ positive
	$k$th powers.
	
	Matolcsi and Ruzsa~\cite[Theorem~2.5 and Corollary~2.6]{MR} proved that
	if $B\subseteq\Z/q\Z$ contains no two elements whose difference is a
	nonzero cube modulo $q$, then $|B|\le q^{1-\kappa}$, where
	$\kappa=0.1195357\ldots$. Our construction, like theirs, uses compatible
	witnesses for cubic residues modulo powers of each prime.
	
	Nin\v{c}evi\'{c} and Slijep\v{c}evi\'{c}~\cite[Theorems~1 and~2]{NS}
	studied the van der Corput property for the values of odd polynomials
	$h\in\Z[x]$ of degree $k\ge3$ with positive leading coefficient.
	They constructed witnesses of degree at most $N$ with nonnegative coefficients
	and $a_0\ll_h(\log N)^{-1/k}$. For $h(n)=n^k$ they also proved that every
	such witness satisfies $a_0\gg_k 1/\log N$.
	
	By allowing the cosine coefficients to have either sign, we obtain a power
	saving for cubes.
	
	\begin{theorem}\label{thm:cubic-witness}
		Let $X\ge1$ and let $\varepsilon>0$. There exists a trigonometric
		polynomial with real coefficients of the form
		\begin{equation}\label{P-witness}
			P_X(\theta)=a_0+\sum_{1\le n\le X}a_n\cos(2\pi n^3\theta)
		\end{equation}
		such that
		\begin{equation}\label{eq:witness-properties}
			P_X(\theta)\ge0\quad(\theta\in\R/\Z),\qquad
			P_X(0)=1,\qquad 0\le a_0\ll_\varepsilon X^{-1/20+\varepsilon}.
		\end{equation}
	\end{theorem}
	
	Taking $X=\lfloor N^{1/3}\rfloor$ and applying
	Theorem~\ref{thm:cubic-witness} with $3\varepsilon$ in place of $\varepsilon$
	gives a witness of degree at most $N$ with
	$a_0\ll_\varepsilon N^{-1/60+\varepsilon}$.
	Applying~\eqref{eq:intro-vdc-bound} gives the following corollary.
	
	\begin{corollary}\label{cor:cube-differences}
		Let $N\ge1$ be an integer and let $A\subseteq[N]$ have no nonzero cube
		difference. Then for every $\varepsilon>0$,
		\[
		|A|\ll_\varepsilon N^{59/60+\varepsilon}.
		\]
	\end{corollary}
	
	Previously, the best upper bound was due to
	Adajar et al.~\cite[Theorem~1.1]{Adajar}. They combined the method of Green
	and Sawhney with exponential sum estimates in the style of Rice~\cite{Rice}
	to show that for every $0<\mu<1/2$ there is a $c_\mu>0$ such that
	\[
	|A|\le N\exp\bigl(-c_\mu(\log N)^\mu\bigr)
	\]
	for all sufficiently large $N$.
	
	In the first version of their
	paper on squares~\cite[Section~1]{FL}, Fan and Lott announced the van der Corput power saving for cube differences that we prove here. After they
	submitted that version to arXiv, OpenAI released a preprint claiming a power saving for intersective polynomial differences by a different
	method~\cite[Theorem~1.1]{OAIpoly}. In the cubic case,
	Corollary~\ref{cor:cube-differences} gives $|A|\ll_c N^{1-c}$ for every
	$0<c<1/60$, whereas the parameters in that preprint give
	$c=c_3\le [72006\,(1201!)^{48001}]^{-1}$.\footnote{Take $k=3$ in~\cite[Sections 2--4 and 7]{OAIpoly}. In the
		notation of that paper, $\delta=1/2400$, $4\delta L>2$, and $t$ is the
		least positive integer with $(t+1)\delta/2>5$. Thus $L\ge1201$, $t=24000$,
		and $d=(L!)^{2t+1}\ge(1201!)^{48001}$. With $\zeta=1/18000$, the proof of
		Proposition~7.2 takes $a_0\le[2(3+2/\zeta)]^{-1}=1/72006$, and the proof of
		Theorem~1.1 takes $c_3=a_0/d$.}
	
	Let $D_3(N)$ denote the largest size of a subset of $[N]$ with no nonzero
	cube difference. Using a construction of Ruzsa~\cite[Theorem~2]{RuzsaSquares}
	and a computer-assisted search, Lewko~\cite[Theorem~4]{Lewko} proved that
	\[
	D_3(N)\gg N^\beta,\qquad
	\beta=\frac23+\frac{\log 14}{3\log 91}=0.861681\ldots.
	\]
	There is still a considerable gap between this exponent $\beta$ and the
	exponent $59/60$ in Corollary~\ref{cor:cube-differences}.
	One may ask whether $D_3(N)\sim C N^\gamma$ for some $C>0$ and
	$\gamma\in[\beta,59/60]$. A weaker conjecture is that the limit
	\[
	\lim_{N\to\infty}\frac{\log D_3(N)}{\log N}
	\]
	exists. Ruzsa~\cite[Section~4]{RuzsaSquares} asked the corresponding
	question for squares.

	\section*{Acknowledgments}
	Andrew Lott and Ákos Magyar thank the Erdős Center of the HUN-REN Alfréd Rényi Institute of Mathematics
	for its support and hospitality during the spring 2026 Simons Semester.
	Máté Matolcsi and Ákos Magyar were supported by grant NKFIH-154121.
	Imre Z.~Ruzsa was supported by grant NKFIH-146387.
	This work was funded by a Simons Foundation International grant
	administered by the Simons Foundation [SFI-MPS-TSM-00026118, Á.M.].
	
	\section*{AI disclosure}
	A large language model (LLM) was used as a sounding board during the developmental phase of the project. The LLM was also used for copyediting. However, the final argument is the authors' own, and the manuscript was written and verified by the authors, who take full responsibility for the contents of the proof.

	\section{Sketch of the proof.}
	
	Our starting point is to define local weight functions $W_{p,K}$ on the groups $\Z/p^K\Z$, where $p$ is prime and $K\in\N$, attached to a parameter $0<\gamma_p<1$. These functions are densities, relative to uniform measure, of geometric probability distributions defined by Pr$(p^j\mid n)=\gamma_p^j$. Extending them as $p^K$-periodic functions, we show that the corresponding exponential sums
	\[
	F_{p,m}(a)=\frac1{p^m}\sum_{n\bmod p^m} W_{p,K}(n)\e(an^3/p^m)\qquad(K=\lceil m/3 \rceil)
	\]
	are nonnegative. With suitable optimized choices of $\gamma_p$, these sums recover the self-compatible modular witness functions constructed by Matolcsi--Ruzsa; see \cite[Lemma 2.4 and Theorem 2.5]{MR} and the remark below. We will use simpler parameter choices. The local weights $W_{p,K}$ are constant on the ``valuation shells" $\{n\in\Z/p^K\Z:v_p(n)=k\}$, where $v_p(n)=j$ if $p^j\mid\mid n$. Thus they can be expressed as linear combinations of the Ramanujan sums $c_{p^k}(n)$:
	\[
	W_{p,K}(n)=\sum_{k=0}^K b_{p^k} c_{p^k}(n),
	\]
	with coefficients $b_{p^k}$ independent of the underlying group $\Z/p^K\Z$. This allows one to define the coefficients $b_d=\prod_{p^k\parallel d}b_{p^k}$ and the global truncated weight function
	\[
	w_R(n) := \sum_{\substack{d\ge1\\b_d\geq R^{-1}}} b_d\,c_d(n).
	\]
	By exploiting the multiplicative properties of the coefficients $b_d$ and the functions $c_d(n)$ we obtain the crucial estimate
	\[
	S_R(a,q):=\mathbb{E}_n w_R(n)\e(an^3/q) \geq -C_\varepsilon R^{-1}q^\varepsilon, \qquad (a,q)=1
	\]
	where by $\mathbb{E}_n f(n)$ we denote the mean value of a periodic function over its period. The reason behind this estimate is that $S_R(a,q)$ can be approximated via the complete exponential sum
	\[
	\E_n\, W_{D}(n)\, \e(an^3/q)=\prod_{p^m\parallel q}F_{p,m}(a_p)\geq0 \qquad (D=\prod_{p^m\parallel q} p^{\lceil m/3\rceil})
	\]
	where $a_p\equiv a(q/p^m)^{-1}\pmod{p^m}$ for $p^m\parallel q$, and
	$W_D(n)=\sum_{d\mid D}b_dc_d(n)=\prod_{p^K\parallel D}W_{p,K}(n)$.
	
	For $X\ge1$ and $R\ge1$, define the weighted exponential sum
	\[
	T_X(\theta)=\sum_{1\le n\le X}h_X(n)w_R(n)e^{2\pi i n^3\theta},
	\qquad
	h_X(t)=\frac{6t^2}{X^3}\left(1-\frac{t^3}{X^3}\right)
	\quad(0\le t\le X).
	\]
	We will take $R=X^{1/20}$. A major/minor arc argument then yields the uniform estimate
	\[
	T_X(0)=1+o(1),\qquad \Re T_X(\theta)\ge-C_\eps X^{-\frac{1}{20}+\eps}.
	\]
	Taking $\eta=C_\varepsilon X^{-1/20+\varepsilon}$ with a sufficiently large constant, we add $\eta$ and normalize:
	\[
	P_X(\theta):=\frac{\Re T_X(\theta)+\eta}{T_X(0)+\eta}.
	\]
	Then $P_X\ge 0$, $P_X(0)=1$, and its constant Fourier coefficient is
	\[
	a_0=\frac{\eta}{T_X(0)+\eta}\ll_\varepsilon X^{-1/20+\varepsilon},
	\]
	providing the witness polynomial in Theorem~\ref{thm:cubic-witness}. A short, well-known argument then gives Corollary~\ref{cor:cube-differences}.
	
	\section{The arithmetic factors} We define the finite exponential sums that appear in the major-arc approximation to $T_X(\theta)$ and prove their required properties. We begin with nonnegative local modular witnesses attached to primes, following the modular construction of Matolcsi--Ruzsa; see \cite[Section 2]{MR}.
	
	\subsection{Local weight functions attached to a prime.} Fix a prime $p$ and a parameter $0<\ga=\ga_p<1$, whose value will be determined later. For $K\ge1$, consider the probability distribution on $\Z/p^K\Z$ that is uniform on each valuation shell $\{v_p(n)=j\}$ , defined by
	\begin{equation}\label{eq:geometric-tails}
		\Pr(p^j\mid n)=\gamma^j\qquad(0\le j\le K).
	\end{equation}
	The uniform distribution corresponds to $\gamma=p^{-1}$. On $\Z/p^K\Z$, we use the valuation: $v_p(n)=j$ when $p^j\parallel n$ and $j<K$, while $v_p(0)=K$. Then
	\begin{equation}\label{local-prob}
		\Pr(v_p(n)=j)=(1-\gamma)\gamma^j\quad(0\le j<K),
		\qquad \Pr(v_p(n)=K)=\gamma^K.
	\end{equation}
	We define the local weight function $W_{p,K}$ on this group as
	the density of this probability distribution relative to uniform distribution, namely as
	\begin{equation}\label{eq:local-density}
		W_{p,K}(n)=
		\begin{cases}
			\displaystyle\frac{p(1-\gamma)}{p-1}(p\gamma)^j,
			&v_p(n)=j<K,\\[5pt]
			(p\gamma)^K,&p^K\mid n.
		\end{cases}
	\end{equation}
	We consider this function as a $p^K$-periodic function on the integers and set $W_{p,0}=1$. It is clear from the construction that it has mean value one over any of its periods, i.e. for every $m\ge K$, the numbers
	$p^{-m}W_{p,K}(n)$, $n\bmod p^m$ also form a probability distribution.
	The probability distributions modulo $p^K$ are compatible as $K$ varies, namely one has the compatibility property
	\[
	\frac1p\sum_{t=0}^{p-1}W_{p,K+1}(n+tp^K)=W_{p,K}(n)
	\qquad(n\bmod p^K).
	\]
	If $v_p(n)=j<K$, all lifts $n+tp^K$ of $n$ have valuation $j$, and their weights equal $W_{p,K}(n)$. For the zero residue modulo $p^K$, one lift has weight $(p\gamma)^{K+1}$ and the other $p-1$ lifts have weight $p(1-\gamma)(p\gamma)^K/(p-1)$. Hence their average is
	\[
	\frac1p\left((p\gamma)^{K+1}+p(1-\gamma)(p\gamma)^K\right)
	=(p\gamma)^K=W_{p,K}(0).
	\]
	This compatibility underlies the construction of the global weight on $\Z$.
	
	\subsection{Positivity of local modular witness functions.}
	
	For $p\nmid a$ and $s\ge1$, define
	\[
	U_{p^s}(a)=\frac1{\varphi(p^s)}
	\sum_{\substack{u\bmod p^s\\p\nmid u}}\e(au^3/p^s),
	\]
	where $\varphi$ is the Euler totient function. The change of variable $u\mapsto-u$ shows that these are
	real and lie in $[-1,1]$. We include a proof of the following standard estimate for completeness.
	
	\begin{lemma}\label{lem:unit}
		If $p\ne3$, then $U_{p^s}(a)=0$ for $s\ge2$. If $p=3$, then
		$U_{3^s}(a)=0$ for $s\ge3$. Moreover, for $p\ne3$,
		\[
		U_p(a)\ge-\frac{1+2\sqrt p}{p-1}.
		\]
	\end{lemma}
	\begin{proof}
		For $p\ne3$ and $s\ge2$, first partition the units $u$ into residue classes modulo $p^{s-1}$. On the residue class represented by $u$, replacing $u$ by $u+tp^{s-1}$ multiplies the phase by $\e(3au^2t/p)$. Summing over $t\bmod p$ gives zero.
		
		For $p=3$ and $s\ge3$, the permutation $u\mapsto u+3^{s-2}$ of
		the units changes the phase by $\e(au^2/3)=\e(a/3)\ne1$.
		The other terms in the binomial expansion are divisible by $3^s$,
		so the sum again vanishes.
		
		Finally, the complete sum $\sum_{x\bmod p}\e(ax^3/p)$ has absolute
		value at most $2\sqrt p$. Indeed, when $p\equiv2\pmod3$, the map $u\mapsto u^3$ is a permutation and the sum is zero. When $p\equiv1\pmod3$, the cubic exponential sum is the sum of the two Gauss sums associated with the nontrivial cubic characters, each of absolute value $\sqrt p$.
		
		Subtracting the term $x=0$ and dividing by $p-1$ proves the asserted lower bound.
	\end{proof}
	
	For $m\ge1$, put $K=\lceil m/3\rceil$ and define the local modular witness functions
	\begin{equation}\label{eq:local-transform}
		F_{p,m}(a)=\frac1{p^m}\sum_{n\bmod p^m}
		W_{p,K}(n)\e(an^3/p^m),\qquad F_{p,0}(a)=1.
	\end{equation}
	
	\begin{lemma}[Local positivity]\label{lem:localpositive}
		For any $0<\gamma<1$ and $p\nmid a$, we have
		\begin{equation}\label{eq:F-recurrence}
			F_{p,m+3}(a)=\gamma F_{p,m}(a)\quad(m\ge0),
		\end{equation}
		with initial values
		\[
		F_{p,1}(a)=\gamma+(1-\gamma)U_p(a),\qquad
		F_{p,2}(a)=\gamma+(1-\gamma)U_{p^2}(a).
		\]
		In particular, with the choices
		\begin{equation}\label{eq:gamma}
			\gamma_3=\frac12,\qquad
			\gamma_p=\frac{1+2\sqrt p}{p+2\sqrt p}\quad(p\ne3),
		\end{equation}
		one has $F_{p,m}(a)\ge0$ for every $m\ge0$ and every $a$ coprime
		to $p$.
	\end{lemma}
	
	\begin{proof}
		Note that if $p^K\mid n$ then $\e(an^3/p^m)=1$.
		If $v_p(n)=j<K$, write $n=p^ju$. The conditional average of the phase on this shell is $U_{p^{m-3j}}(a)$. Thus, summing over $0\le j<K$ gives
		\begin{equation}\label{eq:F-shells}
			F_{p,m}(a)=\gamma^K+
			(1-\gamma)\sum_{j=0}^{K-1}\gamma^j U_{p^{m-3j}}(a).
		\end{equation}
		The first term is the contribution of the residues satisfying $p^K\mid n$.
		
		For $m\ge1$, separating the $j=0$ term for $F_{p,m+3}$ and shifting
		the remaining indices gives
		\[
		F_{p,m+3}(a)=(1-\gamma)U_{p^{m+3}}(a)+\gamma F_{p,m}(a).
		\]
		The term $U_{p^{m+3}}(a)$ vanishes by Lemma~\ref{lem:unit}. The lemma gives $F_{p,3}(a)=\gamma$, proving \eqref{eq:F-recurrence} also
		for $m=0$. Formula~\eqref{eq:F-shells} gives the stated initial
		values for $m=1,2$.
		
		For $p\ne3$, our choice in \eqref{eq:gamma} satisfies
		\[
		\frac{\gamma_p}{1-\gamma_p}=\frac{1+2\sqrt p}{p-1}.
		\]
		Hence Lemma~\ref{lem:unit} implies $F_{p,1}(a)\ge0$, and it is easy to see that $F_{p,2}(a)=\gamma_p>0$.
		
		For $p=3$, the bound $U_3,U_9\ge-1$ and the choice $\gamma_3=1/2$ suffice for both initial values to be nonnegative. Finally, the recurrence \eqref{eq:F-recurrence} implies positivity at every higher modulus $p^m$.
	\end{proof}
	
	\textit{Remark.} The sums $F_{p,m}(a)$ can be evaluated explicitly. If $m=3j+r$ with $r\in\{0,1,2\}$ and $p\nmid a$, then
	\begin{equation}\label{eq:F-explicit}
		F_{p,m}(a)=
		\begin{cases}
			\gamma_p^j,&r=0,\\
			\gamma_p^j\,(\gamma_p+(1-\gamma_p)U_p(a)),&r=1,\\
			\gamma_p^j\,(\gamma_p+(1-\gamma_p)U_{p^2}(a)),&r=2.
		\end{cases}
	\end{equation}
	We have chosen the parameters $\gamma=\ga_p$ above for their simplicity, and have not optimized them. If
	\[
	t_p=\max\bigl(0,\max_{p\nmid a}(-U_p(a)),
	\max_{p\nmid a}(-U_{p^2}(a))\bigr),
	\]
	then $F_{p,m}(a)\ge0$ for all $m\ge0$ and $p\nmid a$ if and only if $\gamma\ge t_p/(1+t_p)$.
	Extend \eqref{eq:local-transform} to all numerators $a$ by the same formula. One may replace $W_{p,\lceil m/3\rceil}$ by $W_{p,m}$ in that formula: the weights agree below valuation $\lceil m/3\rceil$, and the remaining tail has the same mass and phase $1$. Compatibility then gives $F_{p,m}(pa)=F_{p,m-1}(a)$ for $m\ge1$.
	For $p\ne3$, choosing $\gamma=t_p/(1+t_p)$ recovers the self-compatible family in \cite[Lemma 2.4]{MR}, with the base witness averaged under multiplication by nonzero cubes. For $p=3$, the optimal choice in this family is
	$\gamma=\cos(\pi/9)/(1+\cos(\pi/9))$; at modulus $27$ it gives the invariant optimal base witness used in \cite[Theorem 2.5]{MR}, and the higher powers agree with its recursive extension. Thus the geometric probability laws give a realization of those modular witnesses.
	We keep the simpler choices in \eqref{eq:gamma}; optimizing them does not improve the exponent obtained from the estimates below. The properties of the parameters $\ga_p$ that we need are summarized in the following lemma.
	\begin{lemma}\label{lem:gamma-bounds} For all primes $p$, we have
		\begin{equation}\label{eq:gamma-prop}
			(i)\ \  1/p< \ga_p\leq p^{-1/3},\quad\text{and}\quad(ii)\ \ \ga_p \ll p^{-1/2}.
		\end{equation}
	\end{lemma}
	\begin{proof} The second assertion and the inequality $1/p<\gamma_p$ follow immediately from \eqref{eq:gamma}. The upper bound in (i) is also immediate for $p=3$. For $p\ne3$, that upper bound is equivalent to
		\[
		2\sqrt{p}\,(p^{1/3}-1) \leq p-p^{1/3} = (p^{1/3}-1)(p^{2/3}+p^{1/3}),
		\]
		which follows from the arithmetic--geometric mean inequality $2\sqrt p\le p^{2/3}+p^{1/3}$.
	\end{proof}
	
	\subsection{Ramanujan sums and the global weight function \texorpdfstring{$w_R$}{wR}.}
	
	Recall the definition of the Ramanujan sums
	\begin{equation}\label{eq:ramanujan}
		c_d(n)=\sum_{\substack{r\bmod d\\(r,d)=1}}\e(rn/d)
		=\sum_{f\mid d}\mu(d/f)f\1_{f\mid n},\qquad c_1(n)=1.
	\end{equation}
	In particular for $d=p^k$
	\[
	c_{p^k}(n)=p^k\1_{p^k\mid n}-p^{k-1}\1_{p^{k-1}\mid n}
	\qquad(k\ge1).
	\]
	We now express $W_{p,K}$ as a linear combination of $c_{p^k}$, $0\le k\le K$. These $K+1$ functions form an orthogonal basis for the space of functions constant on valuation shells, with respect to the inner product
	\[
	\langle f,g\rangle_K=\frac1{p^K}\sum_{n\bmod p^K}
	f(n)\overline{g(n)}.
	\]
	To be more precise, one has
	\begin{equation}\label{eq:ramanujan-orthogonality}
		\langle c_{p^k},c_{p^\ell}\rangle_K
		=\begin{cases}\varphi(p^k),&k=\ell,\\0,&k\ne\ell,
		\end{cases}
		\qquad 0\le k,\ell\le K.
	\end{equation}
	Indeed, after expanding the two Ramanujan sums,
	a pair of characters has nonzero contribution only if their frequencies
	agree modulo one. Reduced fractions with denominators $p^k$ and
	$p^\ell$ cannot agree mod 1 unless $k=\ell$; in that case there are exactly $\varphi(p^k)$ matching pairs.
	
	Thus one may write
	\begin{equation}\label{eq:W-expansion}
		W_{p,K}(n)=\sum_{k=0}^K b_{p^k}c_{p^k}(n),\qquad
		b_{p^k}=\frac{\langle W_{p,K},c_{p^k}\rangle_K}{\varphi(p^k)}.
	\end{equation}
	The constant coefficient is $b_1=1$, since $W_{p,K}$ has mean one.
	For $1\le k\le K$, direct calculation (or the probability interpretation) gives
	\begin{align*}
		\langle W_{p,K},c_{p^k}\rangle_K
		&=p^k\Pr(p^k\mid n)-p^{k-1}\Pr(p^{k-1}\mid n)\\
		&=p^k\gamma_p^k-p^{k-1}\gamma_p^{k-1}.
	\end{align*}
	Dividing by $\varphi(p^k)=p^{k-1}(p-1)$ yields
	\begin{equation}\label{eq:b}
		b_{p^k}=\frac{p\gamma_p-1}{p-1}\gamma_p^{k-1}
		\qquad(k\ge1).
	\end{equation}
	A crucial fact is that the coefficients $b_{p^k}$ do not depend on the truncation level $K\ge k$. This allows us to extend them multiplicatively, defining $b_d:=\prod_{p^k\parallel d}b_{p^k}$ with $b_1=1$. Note that for $p\ne3$, our choice of the $\ga_p$ parameters gives
	$b_p=2/(\sqrt p+2)$, while $b_3=1/4$.  Since $0<b_p\le\gamma_p$, we obtain
	\begin{equation}\label{eq:b-bound}
		0<b_d\le\prod_{p^k\parallel d}\gamma_p^k\le d^{-1/3}.
	\end{equation}
	
	For $R\ge1$, define the global weight function
	\begin{equation}\label{eq:ramanujan-weight}
		w_R(n):=\sum_{d\in\mathcal E_R}b_d\,c_d(n),\qquad \mathcal E_R:=\{d\ge1:b_d\ge R^{-1}\}.
	\end{equation}
	By \eqref{eq:b-bound}, $\mathcal E_R\subseteq[1,R^3]$, so the defining sum is finite. Although every local weight $W_{p,K}$ is nonnegative, the truncated global weight $w_R$ may have either sign.
	
	\subsection{Lower bounds for the global witness function at rational points \texorpdfstring{$a/q$}{a/q}.}
	
	Our aim is to show that the periodic mean of $w_R(n)\e(an^3/q)$ is bounded below by a small negative quantity. We write $\E_n f(n)$ for the mean of a periodic function $f:\Z\to\mathbb C$ over any of its periods; the value is independent of the period chosen.
	
	\begin{proposition}\label{prop:complete} Let $(a,q)=1$. For every $\varepsilon>0$, uniformly for all reduced fractions $a/q$ and $R\geq 1$, we have
		\begin{equation}\label{eq:w-lowerbound}
			S_R(a,q):=\E_n\,w_R(n)\e(an^3/q)\ge -C_\varepsilon R^{-1}q^\varepsilon.
		\end{equation}
		Moreover $S_R(0,1)=1$.
	\end{proposition}
	Let
	\begin{equation}\label{eq:periodic-means}
		C_d(a,q):=\E_n\,c_d(n)\e(an^3/q)
	\end{equation}
	denote the mean value of the periodic phase function $c_d(n)\e(an^3/q)$ over its period. By \eqref{eq:ramanujan-weight},
	\[
	S_R(a,q)=\sum_{d\in\mathcal E_R}b_dC_d(a,q).
	\]
	The means $C_d(a,q)$ and $S_R(a,q)$ are real: replacing $n$ by $-n$ conjugates the phase and leaves the weights unchanged. The key is to show that $S_R(a,q)$ is a truncated divisor sum.
	
	\begin{lemma}\label{lem:support}
		For $(a,q)=1$, let
		\begin{equation}\label{eq:support}
			D(q):=\prod_{p^m\parallel q}p^{\lceil m/3\rceil},\qquad D(1)=1.
		\end{equation}
		Then
		\[
		C_d(a,q)=0\quad\hbox{unless }d\mid D(q),\]
		and $|C_d(a,q)|\le2^{\omega(d)}$, where $\omega(d)$ is the number of distinct prime factors of $d$.
	\end{lemma}
	\begin{proof}
		A common period of the functions $c_d(n)$ and $\e(an^3/q)$ is $\Delta=[d,q]$. For every prime $p\mid\Delta$,
		write $k_p=v_p(d)$, $m_p=v_p(q)$, and
		$\ell_p=\max(k_p,m_p)$. If $p\mid q$, set
		\[
		q_p=q/p^{m_p},\qquad a_p\equiv a q_p^{-1}\pmod{p^{m_p}}.
		\]
		The inverse $q_p^{-1}\pmod{p^{m_p}}$ exists, and $a_p$ is coprime to $p$. The additive character and the Ramanujan sum factor as
		\[
		\e(an^3/q)=\prod_{p\mid q}\e(a_p n^3/p^{m_p}),\qquad
		c_d(n)=\prod_{p\mid d}c_{p^{k_p}}(n).
		\]
		Indeed, the first identity follows from the fact that
		$1/q\equiv\sum_{p\mid q}q_p^{-1}/p^{m_p}\pmod1$. The
		Chinese remainder theorem applied to $\Delta=\prod_{p\mid\Delta}p^{\ell_p}$ yields the factorization
		\begin{equation}\label{eq:Cd-CRT}
			C_d(a,q)=\prod_{p\mid\Delta}\mathcal L_p(k_p,m_p;a_p),
		\end{equation}
		where, for $m\ge1$ and $\ell=\max(k,m)$,
		\begin{equation}\label{eq:local-Cd-factor}
			\mathcal L_p(k,m;a_p)=\frac1{p^{\ell}}
			\sum_{x\bmod p^{\ell}}c_{p^k}(x)\e(a_p x^3/p^m),
		\end{equation}
		with $c_{p^0}=1$. For $m=0$, the exponential is interpreted as $1$ and no numerator is needed.
		
		If $m=0$ and $k\ge1$, the local factor is the mean of $c_{p^k}$, hence zero. Suppose now that $m\ge1$ and $k\ge1$. For $0\le r\le k$,
		put
		\[
		A_r=\frac{p^r}{p^\ell}
		\sum_{\substack{x\bmod p^\ell\\p^r\mid x}}
		\e(a_p x^3/p^m).
		\]
		This is a normalized exponential sum, so $|A_r|\le1$. The
		identity $c_{p^k}(n)=p^k\1_{p^k\mid n}-p^{k-1}\1_{p^{k-1}\mid n}$ implies $\mathcal L_p(k,m;a_p)=A_k-A_{k-1}$.
		If $k>\lceil m/3\rceil$, then $3(k-1)\ge m$, so the phase is one
		on both supports and $A_k=A_{k-1}=1$. This proves the asserted support property. It also gives $|\mathcal L_p|\le2$ when $k\ge1$;
		when $k=0$, the local factor is an average of numbers of modulus
		one and has absolute value at most one. The stated bound follows from
		\eqref{eq:Cd-CRT}.
	\end{proof}
	
	\begin{proof}[Proof of Proposition~\ref{prop:complete}]
		Let $D=D(q)$. Multiplicativity and \eqref{eq:W-expansion} give
		\[
		\sum_{d\mid D}b_d\, c_d(n)
		=\prod_{p^m\parallel q}W_{p,\lceil m/3\rceil}(n).
		\]
		Thus, by the Chinese remainder theorem,
		\begin{equation}\label{eq:full-positive-product}
			\begin{split}
				\sum_{d\mid D}b_d\,C_d(a,q)
				&=\frac1q\sum_{n\bmod q}
				\left(\prod_{p^m\parallel q}W_{p,\lceil m/3\rceil}(n)\right)
				\e(an^3/q)\\
				&=\prod_{p^m\parallel q}F_{p,m}(a_p)\ge0.
			\end{split}
		\end{equation}
		Indeed, each local weight function $W_{p,\lceil m/3\rceil}$ is $p^m$-periodic and Lemma~\ref{lem:localpositive} applies with $a=a_p$.
		
		By Lemma~\ref{lem:support},
		\[
		S_R(a,q)=\sum_{\substack{d\mid D\\d\in\mathcal E_R}}b_dC_d(a,q).
		\]
		Every omitted coefficient satisfies $b_d<R^{-1}$, so the absolute value of the omitted contribution is at most
		\begin{align*}
			\sum_{\substack{d\mid D\\d\notin\mathcal E_R}}b_d\,|C_d(a,q)|
			&\le R^{-1}\sum_{d\mid D}2^{\omega(d)}\\
			&=R^{-1}\prod_{p^K\parallel D}(1+2K)
			\le R^{-1}\tau(D)^2
			\le R^{-1}\tau(q)^2
			\ll_\varepsilon R^{-1}q^\varepsilon.
		\end{align*}
		Together with the nonnegativity in \eqref{eq:full-positive-product}, this proves \eqref{eq:w-lowerbound}. Finally, for $q=1$ only $d=1$ contributes, so $S_R(0,1)=1$.
	\end{proof}
	
	\section{Lifting to integer cubes via the circle method.}
	
	Using a major/minor arc decomposition, we approximate $\Re T_X(\theta)$ on the major arcs by the product of the real arithmetic factor $S_R(a,q)$ and a nonnegative Archimedean term. The two elementary formulas for $c_d$ in \eqref{eq:ramanujan} serve different purposes: the divisor formula
	is used in the major-arc approximation, and the character formula allows
	us to apply Weyl's inequality on the minor arcs.
	
	\subsection{A preliminary bound.}
	
	The quantity
	\begin{equation}\label{eq:coefficient-cost}
		\mathcal A_R=\sum_{d\in\mathcal E_R}b_d
		\sum_{f\mid d}|\mu(d/f)|f
	\end{equation}
	appears in both the major-arc approximation and the minor-arc estimate. We will use the following bound.
	
	\begin{lemma}\label{lem:coefficient-cost}
		For every $\varepsilon>0$ and $R\ge1$,
		\begin{equation}\label{eq:coeff-bound}
			\mathcal A_R\ll_\varepsilon R^{3+\varepsilon}.
		\end{equation}
	\end{lemma}
	\begin{proof}
		Note that on $\mathcal E_R$ one has $1\le(Rb_d)^{3+\varepsilon}$.
		Extending a sum of non-negative terms to all integers and using
		multiplicativity, we obtain
		\[
		\mathcal A_R
		\le R^{3+\varepsilon}\sum_{d\ge1}b_d^{4+\varepsilon}\,
		\sum_{f\mid d}|\mu(d/f)|f
		=\, R^{3+\varepsilon}\prod_p\left(1+\sum_{k\ge1}
		b_{p^k}^{4+\varepsilon}(p^k+p^{k-1})\right).
		\]
		This Euler product converges. Indeed, put $r_p=p\gamma_p^{4+\varepsilon}$. By Lemma~\ref{lem:gamma-bounds},
		\[
		r_p\le p^{-(1+\varepsilon)/3}<1,
		\]
		and, since $b_{p^k}\le\gamma_p^k$,
		\[
		\sum_{k\ge1}b_{p^k}^{4+\varepsilon}(p^k+p^{k-1})
		\le (1+p^{-1})\sum_{k\ge1}r_p^k
		=(1+p^{-1})\frac{r_p}{1-r_p}.
		\]
		For sufficiently large $p$, the estimate $\gamma_p\ll p^{-1/2}$ gives $r_p=O_\varepsilon(p^{-1-\varepsilon/2})$. Thus the local factor is $1+O_\varepsilon(p^{-1-\varepsilon/2})$, and the finitely many remaining factors are finite. This proves \eqref{eq:coeff-bound}.
	\end{proof}
	
	\subsection{Approximation on the major arcs.} We derive the major-arc approximation by partial summation.
	
	\begin{lemma}\label{lem:major}
		Let $(a,q)=1$ and $\theta=a/q+\beta$. Then
		\[
		\Re T_X(\theta)=S_R(a,q)J(\beta X^3)
		+O\left(\frac{q\,\mathcal A_R}{X}(1+|\beta|X^3)\right),
		\]
		where
		\[
		J(t)=2\int_0^1(1-u)\cos(2\pi tu)\,du
		=\left(\frac{\sin\pi t}{\pi t}\right)^2,
		\qquad J(0)=1.
		\]
		In particular, $0\le J(t)\le1$.
	\end{lemma}
	\begin{proof}
		We use the divisor formula
		\[
		c_d(n)=\sum_{f\mid d}\mu(d/f)f\1_{f\mid n}.
		\]
		For a fixed divisor $f$, the function
		$f\1_{f\mid n}\e(an^3/q)$ has period $fq$ and its mean value over its period is
		\[
		A_f(a,q)=\frac1q\sum_{x\bmod q}\e(af^3x^3/q).
		\]
		Its partial sums, after subtracting the mean times the number of terms, are $O(fq)$. Indeed,
		the sum of the absolute values over a period is $fq$, and
		$|A_f(a,q)|\le1$. After subtracting the mean, complete periods cancel and the remaining terms contribute at most $O(fq)$.
		
		The function $H(t)=h_X(t)\e(\beta t^3)$ vanishes at both endpoints
		of $[0,X]$ and satisfies
		\[
		\int_0^X|H'(t)|\,dt
		\ll X^{-1}+|\beta|X^2
		=X^{-1}(1+|\beta|X^3).
		\]
		Thus summation by parts gives
		\begin{align*}
			\sum_{n\le X}h_X(n)\,f\1_{f\mid n}\e(\theta n^3)
			&=A_f(a,q)\int_0^Xh_X(t)\e(\beta t^3)\,dt\\
			&\quad+O\left(\frac{fq}{X}\,(1+|\beta|X^3)\right).
		\end{align*}
		Multiplying by $b_d\,\mu(d/f)$ and summing over $f\mid d$ and
		$d\in\mathcal E_R$, we obtain the arithmetic factor
		\[
		\sum_{d\in\mathcal E_R}b_d\sum_{f\mid d}\mu(d/f)A_f(a,q)
		=\sum_{d\in\mathcal E_R}b_d\,C_d(a,q)=S_R(a,q).
		\]
		The total error is bounded using \eqref{eq:coefficient-cost}.
		Finally, $S_R(a,q)$ is real. Taking real parts and substituting $u=t^3/X^3$ gives
		\[
		\Re\int_0^Xh_X(t)\e(\beta t^3)\,dt
		=2\int_0^1(1-u)\cos(2\pi\beta X^3u)\,du
		=J(\beta X^3),
		\]
		as required.
	\end{proof}
	
	\subsection{Minor-arc estimates.}
	
	We use the standard cubic Weyl inequality. If $(a,q)=1$ and
	$|\theta-a/q|\le q^{-2}$, then
	\begin{equation}\label{eq:weyl}
		\left|\sum_{n\le Y}\e(\theta n^3+\xi n)\right|
		\ll_\varepsilon Y^{1+\varepsilon}
		(q^{-1}+Y^{-1}+qY^{-3})^{1/4}\qquad(Y\ge1),
	\end{equation}
	uniformly in the real linear coefficient $\xi$, see \cite[Lemma 2.4]{Vaughan} or \cite[equation (1)]{Kerr}.
	
	\begin{lemma}\label{lem:minor}
		Suppose $1\le Q\le X$, $(a,q)=1$, $|\theta-a/q|\le q^{-2}$, and
		$Q<q\le X^3/Q$. Then
		\[
		|T_X(\theta)|
		\ll_\varepsilon \mathcal A_R X^\varepsilon Q^{-1/4}.
		\]
	\end{lemma}
	\begin{proof}
		For every $1\le Y\le X$, the right side of \eqref{eq:weyl}
		is at most
		\[
		C_\varepsilon X^{1+\varepsilon}
		(q^{-1}+X^{-1}+qX^{-3})^{1/4}
		\ll_\varepsilon X^{1+\varepsilon}Q^{-1/4}.
		\]
		Indeed, the expression
		$Y^4/q+Y^3+qY$ increases with $Y$ and the second inequality follows from
		$q>Q$, $X\ge Q$, and $q\le X^3/Q$.
		
		Summation by parts against $h_X$, whose total variation is
		$O(X^{-1})$, gives
		\[
		\left|\sum_{n\le X}h_X(n)\e(\theta n^3+\xi n)\right|
		\ll_\varepsilon X^\varepsilon Q^{-1/4},
		\]
		uniformly in $\xi$. Using the character formula for the Ramanujan sum, we may write
		\[
		T_X(\theta)=\sum_{d\in\mathcal E_R}b_d
		\sum_{\substack{r\bmod d\\(r,d)=1}}
		\sum_{n\le X}h_X(n)\e(\theta n^3+rn/d).
		\]
		Since
		\[
		\varphi(d)=\sum_{f\mid d}\mu(d/f)f\le\sum_{f\mid d}|\mu(d/f)|f,
		\]
		we have $\sum_{d\in\mathcal E_R}b_d\varphi(d)\le\mathcal A_R$. Therefore
		\[
		|T_X(\theta)| \ll_\varepsilon \mathcal{A}_R X^\eps Q^{-1/4}.
		\]
	\end{proof}
	
	\section{Proof of the main result.} We are now in a position to prove the main result. First we show the existence of the witness function $P_X(\theta)$.
	
	\begin{proof}[Proof of Theorem~\ref{thm:cubic-witness}]
		It suffices to treat $0<\varepsilon<1/20$ and sufficiently large $X$. Apply the preceding estimates with $\varepsilon/2$ in place of $\varepsilon$, and put $R=X^{1/20}$, $Q=X^{4/5}$. Given $\theta \in \mathbb{T}=\R/\Z$, Dirichlet approximation with $H=\lfloor X^3/Q\rfloor$ gives a reduced fraction $a/q$ such that
		\[
		1\le q\le H\le X^3/Q,\qquad
		|\beta|:=|\theta-a/q|\le\frac1{qH}
		\le\frac{2Q}{qX^3}.
		\]
		Also $|\beta|\le q^{-2}$ because $q\le H$.
		
		If $q\le Q$, Proposition~\ref{prop:complete} and
		Lemma~\ref{lem:major}, together with $0\le J\le1$, imply
		\[
		\Re T_X(\theta)
		\ge-C_\varepsilon R^{-1}Q^{\varepsilon/2}
		-O(\mathcal A_RQ/X).
		\]
		If $q>Q$, Lemma~\ref{lem:minor} yields
		$|T_X(\theta)|\ll_\varepsilon\mathcal A_RX^{\varepsilon/2}Q^{-1/4}$.
		Lemma~\ref{lem:coefficient-cost} therefore reduces the three
		errors, up to factors bounded by $X^\varepsilon$, to
		\[
		R^{-1}=X^{-1/20},\qquad
		R^3Q/X=X^{-1/20},\qquad
		R^3Q^{-1/4}=X^{-1/20}.
		\]
		This implies the lower bound $\Re T_X(\theta)\ge-C_\varepsilon X^{-1/20+\varepsilon}$ uniformly in $\theta$.
		
		If $\theta=0$, we have  $a=0$, $q=1$, and $\beta=0$ in
		Lemma~\ref{lem:major}. Since $S_R(0,1)=J(0)=1$, it follows
		\[
		T_X(0)=1+O(\mathcal A_R/X)
		=1+O_\varepsilon(X^{-17/20+\varepsilon}).
		\]
		Take $\eta=C_\varepsilon X^{-1/20+\varepsilon}$, with $C_\varepsilon$ large enough that $\Re T_X(\theta)+\eta\ge0$ for all $\theta$. Since $T_X(0)=1+o(1)$, the denominator $T_X(0)+\eta$ is positive for sufficiently large $X$. Define
		\[
		P_X(\theta)=\frac{\Re T_X(\theta)+\eta}{T_X(0)+\eta}.
		\]
		Then $P_X\ge0$, $P_X(0)=1$, and its constant Fourier coefficient is
		\[
		a_0=\frac{\eta}{T_X(0)+\eta}\ll_\varepsilon X^{-1/20+\varepsilon}.
		\]
		All nonzero frequencies are $\pm n^3$ with $1\le n\le X$, and the cosine coefficients are real because $h_X(n)w_R(n)$ is real. For $X$ in a bounded range, use $P_X\equiv1$ and enlarge the implied constant. This completes the proof for all $X\ge1$.
	\end{proof}
	
	Corollary~\ref{cor:cube-differences} now follows from the following short, well-known argument.
	
	\begin{proof}[Proof of Corollary~\ref{cor:cube-differences}]
		Put $X=\lfloor N^{1/3}\rfloor$ and apply Theorem~\ref{thm:cubic-witness} with $3\varepsilon$ in place of $\varepsilon$. Embed $A$ in the cyclic group $G=\Z/L\Z$, where $L=2N+1$, and set $\widehat{\1_A}(j)=\sum_{a\in A}\e(aj/L)$. By orthogonality we have
		\begin{equation}\label{eq:orthogonality}
			\frac1L \sum_{j=0}^{L-1}P_X(j/L)|\widehat{\1_A}(j)|^2=a_0|A|.
		\end{equation}
		Indeed, a nonconstant term could contribute only if
		$a-b\equiv\pm n^3\pmod L$ for some $a,b\in A$. Since
		$|a-b|\le N-1$ and $n^3\le N$, the congruence would be an equality in $\Z$, contradicting our assumption on the set $A$.
		
		The terms on the left of \eqref{eq:orthogonality} are nonnegative; hence, using the fact that $P_X(0)=1$, we have
		\[
		\frac{|A|^2}{L}\le a_0\,|A|.
		\]
		For $A\ne\emptyset$, this implies
		\[
		|A|\le L\,a_0\ll_\varepsilon N X^{-1/20+3\varepsilon}
		\ll_\varepsilon N^{1-1/60+\varepsilon}.
		\]
	\end{proof}

	\bigskip
	\noindent\begin{minipage}{\linewidth}
		\small\raggedright
		\textit{Steve Fan, Andrew Lott, and Ákos Magyar}\\
		Department of Mathematics, University of Georgia, Athens, GA 30602, USA\\
		\texttt{steve.fan@uga.edu}, \texttt{andrew.lott@uga.edu}, \texttt{amagyar@uga.edu}
		
		\medskip
		\textit{Máté Matolcsi}\\
		HUN-REN Alfréd Rényi Institute of Mathematics, Budapest, Hungary\\
		Department of Analysis and Operations Research, Institute of Mathematics,\\
		Budapest University of Technology and Economics, Budapest, Hungary\\
		\texttt{matolcsi.mate@renyi.hu}
		
		\medskip
		\textit{Imre Z.~Ruzsa}\\
		HUN-REN Alfréd Rényi Institute of Mathematics, Budapest, Hungary\\
		\texttt{ruzsa.z.imre@renyi.hu}
	\end{minipage}
	
\end{document}